\documentclass[oneside, reqno]{amsart}

\pdfoutput=1

\usepackage{amsmath, amssymb, amsthm, graphicx}
\usepackage{mathtools}

\graphicspath{{Images/}}

\newtheorem{definition}{Definition}[section]
\newtheorem{theorem}{Theorem}[section]

\newtheorem{lemma}{Lemma}[section]
\newtheorem{proposition}{Proposition}[section]

\newtheorem{corollary}{Corollary}[section]

\begin{document}

\title[Rigidity in de Sitter Space]{A Rigidity Theorem for Spacelike Hypersurfaces in de~Sitter Space}
\author[T. Hasson]{Tristan Hasson}
\address{Tristan Hasson, Department of Mathematical Sciences, University of Durham, Durham DH1 3LE, United Kingdom}
\email{tristan.hasson@durham.ac.uk}
\date{\today}
\begin{abstract}
In this paper we present a rigidity theorem for locally isometric hypersurfaces with a curvature restriction in de~Sitter space. This is an analogue to the case for Riemannian space forms given by Guan and Shen in \cite{GuanShen}.
\end{abstract}
\maketitle

\section{Introduction}
A rigidity theorem for convex regular surfaces in $\mathbb{R}^3$ was proved in 1927 by Cohn-Vossen \cite{CohnVossen}. It states that, given two convex surfaces with the same first fundamental form, these surfaces must differ only by a rigid motion of $\mathbb{R}^3$.

In 2015 Guan and Shen generalised this rigidity theorem to star-shaped hypersurfaces in higher dimensional space forms. Following a similar argument to Guan and Shen, we prove a rigidity result for hypersurfaces in de~Sitter space.

We let $M$ and $\widetilde{M}$ be two spacelike hypersurfaces embedded in de Sitter space (subject to a curvature restriction) for which there exists a local isometry $f:M\rightarrow \widetilde{M}$. We prove that this $f$ must be the restriction of some global isometry of the ambient space. Garding's inequality for hyperbolic polynomials is applied to the second symmetric function of the Weingarten tensor. This together with some integral formulas over the hypersurfaces is used to show that the second fundamental form is preserved under $f$. Hence the first and second fundamental forms are preserved under $f$, which implies $f$ is the restriction of a global isometry.

\section{The Setting}

In this section we introduce the notation and setting for the rigidity theorem. Our ambient space will be de~Sitter space, denoted $(dS^{n+1}, \bar{g})$. We have coordinates $(\rho, \theta)$, where $\theta \in \mathbb{S}^n$. The metric with respect to these coordinates is given by
\begin{equation*}
\bar{g} = ds^2 = -d\rho^2 + \mathrm{cosh}^2(\rho)\sigma,
\end{equation*}
where $\sigma$ denotes the round metric on $\mathbb{S}^n$ with respect to $\theta$. Let $\bar{K}$ denote the sectional curvature of de~Sitter space, and $\bar{\nabla}$ be its Levi-Civita connection. We define the following vector field on de~Sitter space
\begin{equation*}
V=\bar{\phi}(\rho) \frac{\partial}{\partial \rho}.
\end{equation*}
Let $\bar{\phi}(\rho)=\mathrm{cosh}(\rho)$, then define the polar potential function on de~Sitter space as
\begin{equation*}
\bar{\Phi}(\rho) = \int_0^\rho \bar{\phi}(r) dr.
\end{equation*}
Note this means $\bar{\Phi}(\rho)=\mathrm{sinh}(\rho)=\bar{\phi}'(\rho)$.

We will be considering spacelike hypersurfaces in de~Sitter space. These are hypersurfaces whose normal is timelike everywhere, meaning they have positive definite metric everywhere. Spacelike hypersurfaces can be described as graphs over the sphere in de~Sitter space,
\begin{equation*}
\begin{alignedat}{2}
y : \;&\mathbb{S}^n && \to dS^{n+1} \\
&\zeta && \mapsto (y(\zeta), \zeta),
\end{alignedat}
\end{equation*}
with a restriction on the first derivative of $y$ to ensure the normal is timelike. We therefore have canonical coordinates given by $(\zeta)$. In all that follows, we let $(M,g)$ and $(\tilde{M},\tilde{g})$ be locally isometric spacelike hypersurfaces in de~Sitter space, with the isometry denoted
\begin{equation*}
f:M\to\tilde{M}.
\end{equation*}
Let $K$, $\nabla$, $\nu$ and $W(p)$ denote the scalar curvature, Levi-Civita connection, normal and Weingarten map (or shape operator) of $M$ respectively. Similarly define $\tilde{K}$, $\tilde{\nabla}$, $\tilde{\nu}$ and $\tilde{W}(p)$ on $\tilde{M}$. Let $\phi$ and $\tilde{\phi}$ be the restriction of $\bar{\phi}$ to $M$ and $\tilde{M}$ respectively. Similarly let $\Phi$ and $\tilde{\Phi}$ be the restriction of $\bar{\Phi}$.

\section{Hyperbolic Polynomials}

The proof of our rigidity theorem makes vital use of G\r arding's inequality for hyperbolic polynomials. We begin this section by presenting, in our notation, the definitions of G\r arding's theory \cite{Garding}. Results in this brief exposition of G\r arding's theory are given without proof, with references to proofs given. We then go on to apply G\r arding's theory to the case of our rigidity theorem.

\subsection{G\r arding's Theory of Hyperbolic Polynomials}

\begin{definition}
Let $p(x_1, \dots, x_n)$ be a degree $m$ homogeneous polynomial on $\mathbb{R}^n$. If for some $\mathbf{a} \in \mathbb{R}^n$ the univariate polynomial $t \mapsto p(\mathbf{x}+t\mathbf{a})$ has only real solutions for all $\mathbf{x}\in \mathbb{R}^n$, then we say $p$ is hyperbolic with respect to $\mathbf{a}$ (or often simply $p$ is hyperbolic $\mathbf{a}$).
\label{Hyperbolicity}
\end{definition}

\begin{definition}
We define two cones of vectors in $\mathbb{R}^n$:
\begin{equation*}
C(p,\mathbf{a}) \coloneqq \{\mathbf{b}\in \mathbb{R}^n \; | \; \mathrm{the \; roots \; of \;}t \mapsto p(\mathbf{b}+t\mathbf{a}) \mathrm{\; are \; all \; negative}\},
\end{equation*}
\begin{equation*}
S(p, \mathbf{a}) \coloneqq \mathrm{the\ connected\ component\ of\ } \{ \mathbf{x}\in\mathbb{R}^n | p(\mathbf{x})\neq 0\} \mathrm{\ containing\ } \mathbf{a}.
\end{equation*}
\end{definition}

We then have $p$ is hyperbolic with respect to $\mathbf{b}\in C(p,\mathbf{a})=S(p,\mathbf{a})$ (see for example \cite{Renegar}).

\begin{definition}
We define three subspaces of $\mathbb{R}^n$:
\begin{equation*}
E_{C(p,\mathbf{a})}= \{\mathbf{x}\in \mathbb{R}^n \;|\; C(p,\mathbf{a})+\mathbf{x}=C(p,\mathbf{a}) \},
\end{equation*}
\begin{equation*}
L_p = \{ \mathbf{x}\in \mathbb{R}^n \;|\; p(t\mathbf{x}+\mathbf{y})=p(\mathbf{y}),\ \forall\ \mathbf{y},t \},
\end{equation*}
\begin{equation*}
X_{p,\mathbf{a}}= \{ \mathbf{x} \;|\; \mathrm{roots\ of\ } t\mapsto p(t\mathbf{a}+\mathbf{x}) \mathrm{\ are\ all\ zero} \}.
\end{equation*}
\end{definition}

Then we have $E_{C(p,\mathbf{a})}=L_p=X_{p,\mathbf{a}}$ (see for example \cite{Garding}).

\begin{definition}
Let $p(\mathbf{x})$ be a homogeneous polynomial on $\mathbb{R}^n$ of degree $m$. We define the polarized form of $p$ as the unique function $r(\mathbf{x}^1, \dots, \mathbf{x}^m)$, $\mathbf{x}^k \in \mathbb{R}^n$, that is linear in each argument, invariant under permutations of the $\mathbf{x}^k$ and satisfies $r(\mathbf{x}, \dots, \mathbf{x}) = p(\mathbf{x})$.
\label{PolarisedPoly}
\end{definition}

We now give the main result of G\r arding's theory.
\begin{theorem}[G\r arding Inequality for Hyperbolic Polynomials, {\cite[Theorem~5]{Garding}}]
Given $p$ of degree $m>1$, let $p$ be hyperbolic with respect to $\mathbf{a}$ with $p(\mathbf{a})>0$. Let $r$ be the polarised form of $p$, then for $\mathbf{x}^1, \dots, \mathbf{x}^m \in C(p,\mathbf{a})$, the following inequality holds
\begin{equation}
r(\mathbf{x}^1, \dots, \mathbf{x}^m) \geq p(\mathbf{x}^1)^{\frac{1}{m}} \dots p(\mathbf{x}^m)^{\frac{1}{m}}.
\label{GardingEq}
\end{equation}
with equality if and only if $\mathbf{x}^1, \dots , \mathbf{x}^m$ are pairwise proportional modulo $L_p$.
\label{GardingIneq}
\end{theorem}

\subsection{Hyperbolic Polynomials for Rigidity}
\label{HyperbolicRigidity}

This section adapts a result of Guan and Shen applying G\r arding's theory in \cite{GuanShen} to de Sitter space. We will be considering the Weingarten maps (or otherwise shape operators) for all points of a hypersurface, $M$, in de Sitter space, $dS^{n+1}$. These lie in the space of self-adjoint operators (with respect to the induced metric $g$ on $M$) over $T_pM$ for a point $p\in M$, which we will denote $\mathcal{SA}(T_pM)$ (or just $\mathcal{SA}(V)$ for an inner product space $V$).

We define the second symmetric function of a linear map $W$ as
\begin{equation*}
\sigma_2(W)=\sum_{i<j}\kappa_i\kappa_j,
\end{equation*}
where the $\kappa_i$ are the eigenvalues of $W$.
\begin{lemma}
For any given matrix $W$, we have
\begin{equation*}
\sigma_2(W)=\sum_{i<j} w_{ii}w_{jj}-w_{ij}w_{ji}.
\end{equation*}
\label{Sigma2ww}
\end{lemma}

\begin{proof}
First note that the trace of a linear map does not depend on choice of basis. Therefore we have
\begin{equation*}
tr(W)=\sum_i w_{ii} = \sum_i \kappa_{ii}.
\end{equation*}
Then
\begin{equation}
2\sigma_2(W)=2\sum_{i<j}\kappa_i\kappa_j=\sum_{i,j} \kappa_i\kappa_j - \sum_i \kappa_i\kappa_i= tr(W)^2-tr(W^2).
\label{sigma2wproof1}
\end{equation}
Now consider the two traces separately, for components $w_{ij}$ of $W$ under any given basis. We have
\begin{equation}
tr(W)^2=\sum_{i,j} w_{ii} w_{jj} = 2 \sum_{i<j} w_{ii}w_{jj} + \sum_i w_{ii}w_{ii},
\label{sigma2wproof2}
\end{equation}
and
\begin{equation}
tr(W^2)=\sum_{i,j} w_{ij}w_{ji} = 2 \sum_{i<j} w_{ij}w_{ji} + \sum_i w_{ii}w_{ii}.
\label{sigma2wproof3}
\end{equation}
Plugging (\ref{sigma2wproof2}) and (\ref{sigma2wproof3}) into (\ref{sigma2wproof1}) we obtain the result.
\end{proof}

We will now show that $\sigma_2$ is hyperbolic with respect to the identity matrix in $\mathcal{SA}(T_pM)$. We begin by showing the determinant is hyperbolic with respect to the identity in $\mathcal{SA}(T_pM)$ . Since the determinant is of degree $n$, we must show that the univariate polynomial
\begin{equation}
t\mapsto \mathrm{det}(M+tI)
\label{detHyperbolic}
\end{equation}
has $n$ real zeros for all $M\in\mathcal{SA}(T_pM)$. Since we only care if the roots are real, we may as well replace $t$ with $-t$. Then this simply becomes the eigenvalue equation for $M$, but since $M$ is self adjoint we know that all its eigenvalues are real. Therefore (\ref{detHyperbolic}) has $n$ real roots and the determinant is hyperbolic with respect to the identity.

We now state a Lemma of G\r arding \cite{Garding} without proof, followed by a corollary.

\begin{lemma}
Let $p$ be a degree $m$ polynomial on $\mathbb{R}^n$ that is hyperbolic with respect to $\mathbf{a}$. The polynomial
\begin{equation}
q(x_1, \dots, x_n)=\sum_{k=1}^n a_k \frac{\partial}{\partial x_k} p(x_1, \dots, x_n)
\end{equation}
is hyperbolic with respect to $\mathbf{a}$.
\label{GardingLemma1}
\end{lemma}

\begin{corollary}
Given a degree $m>1$ polynomial $p$, hyperbolic with respect to $\mathbf{a}\in \mathbb{R}^n$, the polynomials $\{p_i\}_{i=0}^m$, defined by
\begin{equation*}
p(s\mathbf{a}+\mathbf{x})=\sum_{i=0}^m s^i p_i(\mathbf{x}),
\end{equation*}
are also hyperbolic with respect to $\mathbf{a}$.
\label{SymmetricCor}
\end{corollary}

\begin{proof}
Consider evaluating this polynomial $p(s\mathbf{a}+\mathbf{x})$ at $s=0$. We then have
\begin{equation*}
p(\mathbf{x})=p_0(\mathbf{x})
\end{equation*}
and since $p$ is hyperbolic $\mathbf{a}$, $p_0$ must also be hyperbolic $\mathbf{a}$. Now consider the polynomial
\begin{equation*}
q(s\mathbf{a}+\mathbf{x})=\frac{d}{ds}p(s\mathbf{a}+\mathbf{x})=\sum_{i=1}^m i s^{i-1} p_i(\mathbf{x}).
\end{equation*}
Evaluating this at $s=0$, we obtain
\begin{equation*}
q(\mathbf{x})=p_1(\mathbf{x}).
\end{equation*}
Note that this $q$ is exactly the same $q$ as defined in Lemma $\ref{GardingLemma1}$, so since $p$ is hyperbolic $\mathbf{a}$, we have $q=p_1$ is hyperbolic $\mathbf{a}$. Repeating this process shows that all the $p_i$ are hyperbolic with respect to $\mathbf{a}$.
\end{proof}

Now note that the characteristic polynomial for a matrix satisfies
\begin{equation*}
\begin{split}
(-1)^n \; \mathrm{det}(M-tI)&= \sum_{i=0}^n (-1)^i \; t^{n-i} \; (\mathrm{sum\;of\;all\;combinations\;of\;}i\mathrm{\;eigenvalues}) \\
&= \sum_{i=0}^n (-1)^i \; t^{n-i} \sigma_i(M),
\end{split}
\end{equation*}
where $\sigma_i(M)$ is the $i^{\mathrm{th}}$ symmetric function of the eigenvalues of $M$. Now since the determinant is hyperbolic with respect to the identity, we have $(-1)^n det (M)$ is hyperbolic with respect to $-I$. Then by Corollary \ref{SymmetricCor}, we have that the functions $(-1)^i \; \sigma_i$ are hyperbolic with respect to $-I$. So finally, the symmetric functions, in particular $\sigma_2$, are hyperbolic with respect to the identity.

The following lemma (see \cite{GuanShen}) can easily be derived from the Gauss equation (see for example \cite{DoCarmo}).
\begin{lemma}
Let $K$ and $\bar{K}$ denote the scalar curvatures of $M^n$ and $dS^{n+1}$ respectively. We have yet another expression for $\sigma_2$:
\begin{equation*}
\sigma_2(W)=\frac{n(n-1)}{2}(K-\bar{K}).
\end{equation*}
As a corollary, since the scalar curvature is invariant under local isometries, we have $\sigma_2(W)=\sigma_2(\tilde{W})$.
\label{sigma2KK}
\end{lemma}

From Lemma \ref{sigma2KK} note that if $K>\bar{K}$ for all points on the hypersurface, then $\sigma_2(W)$ is positive for all points on the hypersurface. We now claim that this implies that the Weingarten maps $W(p)$ for every point $p\in M$ all lie in the same hyperbolicity cone. For any point $M \in \mathcal{SA}(T_pM)$, consider the affine line through $M$ in the direction of $I$. Since $\sigma_2$ is hyperbolic with respect to the identity, this affine line must cross the hypersurface $\sigma_2=0$ exactly twice (with multiplicities). Let $t_1< t_2$ be the roots of $t\mapsto \sigma_2(M+tI)$. Since $\sigma_2$ is of degree two, if $t_1\neq t_2$ then $\sigma_2(M+tI)$ changes sign at $t_1$ and at $t_2$, and if $t_1=t_2$ then $\sigma_2(M+tI)$ has the same sign either side of the root. Importantly, note that in either case we have points $\sigma_2(M+tI)$ has the same sign for $t<t_1$ and $t>t_2$, and values $t_1<t<t_2$ have the opposite sign. Now by definition, the points $M+tI$ with $t>t_2$ are in $C(\sigma_2, I)$ and the points $M+tI$ with $t<t_1$ are in $C(\sigma_2, -I)$. Since $C(\sigma_2, I)=S(\sigma_2,I)$, for any affine line $M+tI$, we have $\mathrm{sgn}(\sigma_2(M+tI)) = \mathrm{sgn}(\sigma_2(I)) > 0$ for all points with $t<t_1$ or $t>t_2$, and for any other $t_1<t<t_2$ we have $\mathrm{sgn}(\sigma_2(M+tI))<0$. This means $\sigma_2$ is hyperbolic with respect to any $W\in \mathcal{SA}(T_pM)$ with $\sigma_2(W)>0$. Since we stipulate that $K>\bar{K}$, by Lemma \ref{sigma2KK} $\sigma_2$ is hyperbolic with respect to $W(p)$ for all points $p\in M$. Furthermore, since $\sigma_2(W(p))$ is strictly greater than zero, it always stays in the same connected component of $\{M\in \mathcal{SA}(T_pM) \; | \; \sigma_2(M)\neq 0)\}$. So either $\{W(p) \; | \; p \in M \} \subset S(\sigma_2, I)$ or $\{W(p) \; | \; p \in M \} \subset S(\sigma_2, -I)$.

We will be considering two isometric hypersurfaces $M$ and $\tilde{M}$ and their respective Weingarten tensors $W$ and $\tilde{W}$. We would like Theorem \ref{GardingIneq} to hold for $W(p)$ and $\tilde{W}(f(p))$ for all points $p$ and $f(p)$ identified under the isometry. In order for this to hold, we need that $W(p)$ and $\tilde{W}(f(p))$ both lie in the same hyperbolicity cone, i.e. either $C(\sigma_2, I)$ or $C(\sigma_2, -I)$. Note that $S(\sigma_2, I) = -S(\sigma_2, -I)$.

We now show that by using a global isometry of de Sitter space, we can make sure all Weingarten maps for both hypersurfaces lie in $C(\sigma_2, I)$. Since the hypersurface is spacelike, it must be compact. Therefore, we can find a point $p$ such that the $\rho$ coordinate of $p$ is greater than any other point on the hypersurface. As the hypersurface is described as the graph of a function $f:\mathbb{S}^n\rightarrow \mathbb{R}$, we can calculate the second fundamental form at this point as
\begin{equation*}
B_{ij}=\Gamma^0_{ij}+\frac{\partial^2 f}{\partial \zeta_i \partial \zeta_j},
\end{equation*}
where $\Gamma^0_{ij}$ are Christoffel symbols on de Sitter space running over the $\{\zeta_i\}_{i=1}^n$ coordinates. We have used the fact that $f$ has a global maximum at $p$, meaning that all first derivatives are zero. To find the Weingarten map, we simply contract with the induced inverse metric, $g^{ij}$, of the hypersurface at this point. Calculating the Christoffel symbols, we find
\begin{equation}
\sum_k g^{ik} \Gamma^0_{kj} = \mathrm{cosh}(\rho)\mathrm{sinh}(\rho) \; \delta^i_j,
\label{Gamma0ij}
\end{equation}
where $\delta^i_j$ is the Kronecker delta symbol.

Now consider reflecting this hypersurface in the ``equator'' of de Sitter space, $\rho=0$. This simply amounts to replacing the $\rho$ coordinate of any point $x$ in the hypersurface by $-\rho$. Plugging $-\rho$ in (\ref{Gamma0ij}) we obtain
\begin{equation*}
\mathrm{cosh}(-\rho)\mathrm{sinh}(-\rho) \; \delta^i_j = -\mathrm{cosh}(\rho)\mathrm{sinh}(\rho) \; \delta^i_j.
\end{equation*}
As the $\rho$ coordinate of a point on the hypersurface is given by the function $f$, this reflection is obtained by replacing $f$ with $-f$. So since
\begin{equation*}
\frac{\partial^2 (-f)}{\partial \zeta_i \partial \zeta_j}=-\frac{\partial^2 (f)}{\partial \zeta_i \partial \zeta_j},
\end{equation*}
we have the Weingarten map of the maximal point $p$, under the reflection, is given by
\begin{equation*}
\hat{W}_{ij}=\sum_k g^{ik} \hat{B}_{kj} = -\mathrm{cosh}(\rho)\mathrm{sinh}(\rho) \; \delta^i_j - \sum_k g^{ik} \frac{\partial^2 (f)}{\partial \zeta_k \partial \zeta_j} = -W_{ij}
\label{W=-W}.
\end{equation*}

This shows that by reflecting about the ``equator'' we can move the Weingarten map of the maximum point from $C(\sigma_2, -I)$ to $C(\sigma_2,I)$. So up to global isometry, $W(p), \tilde{W}(f(p)) \in C(\sigma_2, I)$ for all points $p\in M$. Therefore the G\r arding inequality for hyperbolic polynomials holds for all points on the hypersurfaces.

\subsubsection{The Equality Result}
Since $W(p),\tilde{W}(f(p)) \in C(\sigma_2,I)$ for all points $p\in M$, and further $\sigma_2(W)=\sigma_2(\tilde{W})$, we have the inequality
\begin{equation*}
\sigma_{1,1}(W,\tilde{W})\geq \sigma_2(W),
\end{equation*}
where $\sigma_{1,1}(W,\tilde{W})$ denotes the polarised form of $\sigma_2(W)$. As $\sigma_2$ is degree two its polarised form has the following simple form (see \cite{GuanShen}).

\begin{proposition}
The polarised form of $\sigma_2(W)$ as defined in Definition \ref{PolarisedPoly} can be expressed as
\begin{equation*}
\sigma_{1,1}(W,\tilde{W})= \frac{1}{2} \sum_{i,j} \frac{\partial \sigma_2 (W)}{\partial w_{ij}} \tilde{w}_{ij}.
\end{equation*}
\label{Sigma11}
\end{proposition}

Now by Theorem \ref{GardingIneq} we have if
\begin{equation*}
\sigma_{1,1}(W,\tilde{W})= \sigma_2(W)
\end{equation*}
then $W$ and $\tilde{W}$ must be proportional modulo $L_{\sigma_2}$. As such, we would like to know what $L_{\sigma_2}$ actually is, and it turn out that in fact it is trivial.
\begin{proposition}
$L_{\sigma_2}=\mathbf{0}$
\end{proposition}
\begin{proof}
As seen in G\r arding's paper \cite{Garding}, $L_{\sigma_2}$ is the set of $W$ such that the roots of the polynomial
\begin{equation*}
t\mapsto \sigma_2(tI+W)
\end{equation*}
are all zero. Expanding this out, we have
\begin{equation*}
\begin{split}
\sigma_2(tI+W)&=\sum_{i<j} (t+w_{ii})(t+w_{jj})-w_{ij}w_{ji} \\
&=\sum_{i<j} t^2 +(w_{ii}+w_{jj})t+ w_{ii}w_{jj}-w_{ij}w_{ji}.
\end{split}
\end{equation*}
For all the roots of this to be zero, we need
\begin{equation*}
\sum_{i<j} (w_{ii}+w_{jj}) = (n-1) \sum_i w_{ii} =0
\end{equation*}
and
\begin{equation*}
\sum_{i<j} w_{ii}w_{jj} - w_{ij}w_{ji}=0.
\end{equation*}
Put simply, this means we need $\sigma_1(W)=\sigma_2(W)=0$. Now using the definition of $\sigma_2(W)$ in terms of the eigenvalues $\kappa_i$ of $W$, we have
\begin{equation*}
\begin{split}
\sigma_2(W)&=\sum_{i<j} \kappa_i \kappa_j \\
&=\frac{1}{2} \sum_{i\neq j} \kappa_i \kappa_j \\
&=\frac{1}{2} \Big( \sum_{i,j} \kappa_i \kappa_j - \sum_i \kappa_i \kappa_i \Big) \\
&=\frac{1}{2} \Big( \sum_{i} \kappa_i \Big( \sum_j \kappa_j \Big) - \sum_i \kappa_i \kappa_i \Big).
\end{split}
\end{equation*}
Now since $\sigma_1(W)=\sum_i \kappa_i =0$, we have
\begin{equation*}
\sigma_2(W)=-\frac{1}{2} \sum_i \kappa_i \kappa_i.
\end{equation*}
Then clearly the only way for $\sigma_2(W)=0$ is for all the $\kappa_i=0$. Hence $L_{\sigma_2}=\mathbf{0}$.
\end{proof}

Since $\sigma_2(W)=\sigma_2(\tilde{W})$, if $W$ and $\tilde{W}$ are proportional then we must have $W=\tilde{W}$. Hence we have
\begin{equation*}
\sigma_{1,1}(W,\tilde{W})=\sigma_2(W)
\end{equation*}
if and only if $W=\tilde{W}$.

\section{Some Important Integrals}

\subsection{A Few Lemmas}
The following lemmas (see \cite{GuanShen}) give the integrands for the integral equations which follow, these will be vital in the main rigidity proof.

\begin{lemma}
\begin{equation*}
\langle \bar{\nabla}_{e_i} V, e_j \rangle_{dS} + \langle \bar{\nabla}_{e_j} V, e_i \rangle_{dS} = 2 \phi' \bar{g}_{ij}
\end{equation*}
\label{DerivV}
\end{lemma}

\begin{proof}
We will prove this in two stages. First, we show that $(\mathcal{L}_V \bar{g})(e_i,e_j) = \langle \bar{\nabla}_{e_i} V, e_j \rangle_{dS} + \langle \bar{\nabla}_{e_j} V, e_i \rangle_{dS}$, and then that $(\mathcal{L}_V \bar{g})(e_i,e_j) = 2\phi'(\rho)\bar{g}_{ij}$.

So first, since the Lie derivative obeys Leibniz's rule, we have
\begin{equation*}
\mathcal{L}_V (\bar{g}(e_i, e_j)) = (\mathcal{L}_V \bar{g})(e_i,e_j) + \bar{g}(\mathcal{L}_V e_i, e_j) + \bar{g}(e_i, \mathcal{L}_V e_j).
\end{equation*}
Using that the derivative of a function is the directional derivative and that $\mathcal{L}_X Y = [X,Y]$, we have
\begin{equation*}
(\mathcal{L}_V \bar{g})(e_i,e_j) = \langle \mathrm{grad}(\langle e_i, e_j \rangle_{dS}), V \rangle_{dS} - \langle [V, e_i], e_j \rangle_{dS} - \langle e_i, [V, e_j] \rangle_{dS},
\end{equation*}
and taking $V$ into the inner product using compatibiliity with the metric, we obtain
\begin{equation*}
(\mathcal{L}_V \bar{g})(e_i,e_j) = \langle \bar{\nabla}_V e_i, e_j \rangle_{dS} + \langle e_i, \bar{\nabla}_V e_j \rangle_{dS} - \langle [V, e_i], e_j \rangle_{dS} - \langle e_i, [V, e_j] \rangle_{dS}.
\end{equation*}
Now using the symmetry of the connection $\nabla_X Y - \nabla_Y X = [X,Y]$, we have
\begin{equation*}
(\mathcal{L}_V \bar{g})(e_i,e_j) = \langle \bar{\nabla}_{e_i} V, e_j \rangle_{dS} + \langle \bar{\nabla}_{e_j} V, e_i \rangle_{dS}.
\label{VLemma1.5}
\end{equation*}

The second statement $(\mathcal{L}_V \bar{g})(e_i,e_j) = 2\phi'(\rho)\bar{g}_{ij}$, is proved simply by calculation using the Cartan formula for the Lie derivative of a differential form. Hence the lemma is proved.
\end{proof}

\begin{lemma}
Let $\nu$ be the normal to $M$, then we have
\begin{equation*}
\mathrm{Hess}^{M,g}(\Phi)(e_i,e_j) = \phi' g_{ij} + h_{ij} \langle V, \nu \rangle_{dS}.
\end{equation*}
\label{preIntegral}
\end{lemma}

\begin{proof}
We start with the definitions of the Hessian for $\bar{\Phi}$ in de Sitter space and for $\Phi$ in $M$,
\begin{subequations}
\begin{equation}
\mathrm{Hess}^{dS,\bar{g}}(\bar{\Phi})(e_i,e_j) = \langle \mathrm{grad}(\langle \mathrm{grad}(\bar{\Phi}), e_j \rangle_{dS} ), e_i \rangle_{dS} - \langle \mathrm{grad}(\bar{\Phi}), \bar{\nabla}_{e_i} e_j \rangle_{dS},
\label{HessdS}
\end{equation}
\begin{equation}
\mathrm{Hess}^{M,g}(\Phi)(e_i,e_j) = \langle \mathrm{grad}(\langle \mathrm{grad}(\Phi), e_j \rangle_M ), e_i \rangle_M - \langle \mathrm{grad}(\Phi), \nabla_{e_i} e_j \rangle_M.
\label{HessM}
\end{equation}
\end{subequations}
Subtracting (\ref{HessdS}) from (\ref{HessM}), we obtain
\begin{multline*}
\mathrm{Hess}^{M,g}(\Phi)(e_i,e_j) - \mathrm{Hess}^{dS,\bar{g}}(\bar{\Phi})(e_i,e_j) \\
= \langle \mathrm{grad}(\langle \mathrm{grad}(\Phi), e_j \rangle_M ), e_i \rangle_M - \langle \mathrm{grad}(\langle \mathrm{grad}(\bar{\Phi}), e_j \rangle_{dS} ), e_i \rangle_{dS} \\
- \langle \mathrm{grad}(\Phi), \nabla_{e_i} e_j \rangle_M + \langle \mathrm{grad}(\bar{\Phi}), \bar{\nabla}_{e_i} e_j \rangle_{dS}.
\end{multline*}
Now clearly $\langle \mathrm{grad}(\langle \mathrm{grad}(\Phi), e_j \rangle_M ), e_i \rangle_M = \langle \mathrm{grad}(\langle \mathrm{grad}(\bar{\Phi}), e_j \rangle_{dS} ), e_i \rangle_{dS}$, and also $\langle \mathrm{grad}(\Phi), \nabla_{e_i} e_j \rangle_M = \langle \mathrm{grad}(\bar{\Phi}), \nabla_{e_i} e_j \rangle_{dS}$, so we have
\begin{equation*}
\mathrm{Hess}^{M,g}(\Phi)(e_i,e_j) = \mathrm{Hess}^{dS,\bar{g}}(\bar{\Phi})(e_i,e_j) + \langle \mathrm{grad}(\bar{\Phi}), \bar{\nabla}_{e_i} e_j - \nabla_{e_i} e_j \rangle_{dS}.
\end{equation*}
Note that on de Sitter space, $\mathrm{grad}{\bar{\Phi}} = V$, so then from the definition of the Hessian, we obtain
\begin{equation}
\mathrm{Hess}^{M,g}(\Phi)(e_i,e_j) = \langle \bar{\nabla}_{e_i} V, e_j \rangle_{dS} + \langle V, \bar{\nabla}_{e_i} e_j - \nabla_{e_i} e_j \rangle_{dS}.
\label{HessM2}
\end{equation}
Note $\langle V, \bar{\nabla}_{e_i} e_j - \nabla_{e_i} e_j \rangle_{dS} = h_{ij} \langle V, \nu \rangle_{dS}$. Now take this equation and consider swapping $e_i$ and $e_j$ in each term
\begin{equation*}
\mathrm{Hess}^{M,g}(\Phi)(e_j,e_i) = \langle \bar{\nabla}_{e_j} V, e_i \rangle_{dS} + h_{ji} \langle V, \nu \rangle_{dS}.
\end{equation*}
Adding this to (\ref{HessM2}), and using that both the Hessian and second fundamental form are symmetric in $i$ and $j$, we obtain
\begin{equation*}
2 \mathrm{Hess}^{M,g}(\Phi)(e_i,e_j) = \langle \bar{\nabla}_{e_i} V, e_j \rangle_{dS} + \langle \bar{\nabla}_{e_j} V, e_i \rangle_{dS} + 2 h_{ij} \langle V, \nu \rangle_{dS}.
\end{equation*}
Finally by Lemma \ref{DerivV} we obtain the result
\begin{equation*}
\mathrm{Hess}^{M,g}(\Phi)(e_i,e_j) = \phi' g_{ij} + h_{ij} \langle V, \nu \rangle_{dS}.
\end{equation*}
\end{proof}

\subsection{The Integral Equations}
By identifying points of $M$ and $\tilde{M}$ under the isometry $f$, we can treat $\tilde{\phi}$ and $\tilde{\Phi}$ as functions on $M$, i.e. for $x\in M$, $\tilde{\phi}_M: x \mapsto \tilde{\phi}(f(x))$. We will abuse notation and use $\tilde{\phi}$ to denote $\tilde{\phi}_M$ in integrals over $M$. The following theorem is adapted from a lemma in \cite{GuanShen}.

\begin{theorem}
Given an orthonormal frame $\{e_i\}$ on $M$, which can be identified as an orthonormal frame on $\tilde{M}$ under the isometry, the following integral equations hold:
\begin{subequations}
\begin{equation}
\int_M \sum_{i,j} \frac{\partial \sigma_2}{\partial w_{ij}}(W) \; \tilde{\phi}' \; \mathrm{Hess}^{M,g}(\Phi)(e_i, e_j) = \int_M (n-1) \tilde{\phi}' \phi' \sigma_1(W) - 2 \tilde{\phi}' \sigma_2(W) \langle V, \nu \rangle,
\label{Integral0a}
\end{equation}
\begin{equation}
\int_M \sum_{i,j} \frac{\partial \sigma_2}{\partial w_{ij}}(\tilde{W}) \; \tilde{\phi}' \; \mathrm{Hess}^{M,g}(\Phi)(e_i, e_j) = \int_M (n-1) \tilde{\phi}' \phi' \sigma_1(\tilde{W}) - 2 \tilde{\phi}' \sigma_{1,1}(W,\tilde{W}) \langle V, \nu \rangle,
\label{Integral0b}
\end{equation}
\begin{equation}
\int_M \sum_{i,j} \frac{\partial \sigma_2}{\partial w_{ij}}(W) \; \phi' \; \mathrm{Hess}^{M,g}(\tilde{\Phi})(e_i, e_j) = \int_M (n-1) \phi' \tilde{\phi}' \sigma_1(W) - 2 \phi' \sigma_{1,1}(\tilde{W}\!,W) \langle \tilde{V}\!, \tilde{\nu} \rangle,
\label{Integral0c}
\end{equation}
\begin{equation}
\int_M \sum_{i,j} \frac{\partial \sigma_2}{\partial w_{ij}}(\tilde{W}) \; \phi' \; \mathrm{Hess}^{M,g}(\tilde{\Phi})(e_i, e_j) = \int_M (n-1) \phi' \tilde{\phi}' \sigma_1(\tilde{W}) - 2 \phi' \sigma_2(\tilde{W}) \langle \tilde{V}, \tilde{\nu} \rangle.
\label{Integral0d}
\end{equation}
\end{subequations}
\label{IntegralEqn}
\end{theorem}

\begin{proof}
Starting with Lemma \ref{preIntegral} applied to $M$ and $\tilde{M}$, we have
\begin{subequations}
\begin{equation}
\mathrm{Hess}^{M,g}(\Phi)(e_i,e_j) = \phi' g_{ij} + h_{ij} \langle V, \nu \rangle_{dS}
\label{Integral1a}
\end{equation}
\begin{equation}
\mathrm{Hess}^{\tilde{M},\tilde{g}}(\tilde{\Phi})(e_i,e_j) = \tilde{\phi}' \tilde{g}_{ij} + \tilde{h}_{ij} \langle \tilde{V}, \tilde{\nu} \rangle_{dS}
\label{Integral1b}
\end{equation}
\end{subequations}
It is not hard to convince oneself that $\mathrm{Hess}^{\tilde{M},\tilde{g}}(\tilde{\Phi})(e_i,e_j)=\mathrm{Hess}^{M,g}(\tilde{\Phi})(e_i,e_j)$. Now we multiply (\ref{Integral1a}) by $\tilde{\phi}'$ and (\ref{Integral1b}) by $\phi'$. Now take the four combinations of multiplying these two equations by either $\frac{\partial\sigma_2}{\partial w_{ij}}(W)$ or $\frac{\partial\sigma_2}{\partial w_{ij}}(\tilde{W})$, sum over the $i$ and $j$ and integrate over $M$ to obtain
\begin{subequations}
\begin{equation}
\int_M \sum_{i,j} \frac{\partial \sigma_2}{\partial w_{ij}}(W) \; \tilde{\phi}' \; \mathrm{Hess}^{M,g}(\Phi)(e_i,e_j) = \int_M \sum_{i,j} \frac{\partial \sigma_2}{\partial w_{ij}}(W) \Big( \tilde{\phi}' \phi' g_{ij} + \tilde{\phi}' h_{ij} \langle V, \nu \rangle_{dS}\Big),
\label{Integral2a}
\end{equation}
\begin{equation}
\int_M \sum_{i,j} \frac{\partial \sigma_2}{\partial w_{ij}}(\tilde{W}) \; \tilde{\phi}' \; \mathrm{Hess}^{M,g}(\Phi)(e_i,e_j) = \int_M \sum_{i,j} \frac{\partial \sigma_2}{\partial w_{ij}}(\tilde{W}) \Big(\tilde{\phi}' \phi' g_{ij} + \tilde{\phi}' h_{ij} \langle V, \nu \rangle_{dS}\Big),
\label{Integral2b}
\end{equation}
\begin{equation}
\int_M \sum_{i,j} \frac{\partial \sigma_2}{\partial w_{ij}}(W) \; \phi' \; \mathrm{Hess}^{M,g}(\tilde{\Phi})(e_i,e_j) = \int_M \sum_{i,j} \frac{\partial \sigma_2}{\partial w_{ij}}(W) \Big(\phi' \tilde{\phi}' \tilde{g}_{ij} + \phi' \tilde{h}_{ij} \langle \tilde{V}, \tilde{\nu} \rangle_{dS}\Big),
\label{Integral2c}
\end{equation}
\begin{equation}
\int_M \sum_{i,j} \frac{\partial \sigma_2}{\partial w_{ij}}(\tilde{W}) \; \phi' \; \mathrm{Hess}^{M,g}(\tilde{\Phi})(e_i,e_j) = \int_M \sum_{i,j} \frac{\partial \sigma_2}{\partial w_{ij}}(\tilde{W}) \Big(\phi' \tilde{\phi}' \tilde{g}_{ij} + \phi' \tilde{h}_{ij} \langle \tilde{V}, \tilde{\nu} \rangle_{dS}\Big).
\label{Integral2d}
\end{equation}
\end{subequations}
We now focus on the right hand side of (\ref{Integral2a}). Since $\{e_i\}$ is an orthonormal basis, we have $g_{ij}=\delta_{ij}$, so for the first term we have
\begin{equation}
\sum_{i=j} \frac{\partial \sigma_2}{\partial w_{ij}}(W) \tilde{\phi}' \phi'.
\label{Integral2.5}
\end{equation}
The expression for $\frac{\partial \sigma_2}{\partial w_{ij}}(W)$ can easily shown by calculation to be
\begin{equation}
\frac{\partial \sigma_2 (W)}{\partial w_{ij}}=
\begin{cases}
\sum\limits_{k\neq i} w_{kk} & \mathrm{for} \; i=j \\
-w_{ji} & \mathrm{for} \; i\neq j.
\end{cases}
\label{Sigma2Deriv}
\end{equation}
Subsituting this into (\ref{Integral2.5}) gives
\begin{equation*}
\sum_{i} \sum\limits_{k\neq i} w_{kk} \tilde{\phi}' \phi',
\end{equation*}
which is clearly equal to
\begin{equation}
(n-1) \tilde{\phi}' \phi' \sigma_1(W).
\label{Integral3}
\end{equation}
For the second term of the right hand side of (\ref{Integral2a}), note that since $g_{ij}=\delta_{ij}$, we have $h_{ij}=w_{ij}$. Therefore we have
\begin{equation*}
\tilde{\phi}' \langle V, \nu \rangle \sum_{i,j} \frac{\partial \sigma_2(W)}{\partial w_{ij}} w_{ij}.
\end{equation*}
From the expression for $\sigma_{1,1}(W, \tilde{W})$ in Proposition \ref{Sigma11}, this is equal to
\begin{equation*}
2\tilde{\phi}' \langle V, \nu \rangle \sigma_{1,1}(W,W),
\end{equation*}
but as the polarised form of $\sigma_2$, we have $\sigma_{1,1}(W,W)=\sigma_2(W)$. This together with (\ref{Integral3}) gives the right hand side of (\ref{Integral0a}). The remaining three equations are proved along the same lines, hence the theorem is proved.
\end{proof}

\subsection{Symmetry in Tilde}

We have the following proposition from \cite{GuanShen}.

\begin{proposition}
Assume that $W$ is a Codazzi tensor and $\{e_1, \dots, e_n\}$ is an orthonormal frame on $M$, then we have the following identity
\begin{equation*}
\sum_i \Big\langle \mathrm{grad}\Big(\frac{\partial \sigma_2}{w_{ij}}(W)\Big), e_i \Big\rangle=0.
\end{equation*}
\label{gradSigma2}
\end{proposition}

We can now prove our symmetry theorem.

\begin{theorem}
The integral
\begin{equation*}
\int_M \sum_{i,j} \frac{\partial \sigma_2}{\partial w_{ij}}(W) \; \tilde{\phi}' \; \mathrm{Hess}^{M,g}(\Phi)(e_i, e_j)
\label{Sym1}
\end{equation*}
is invariant under switching $\tilde{\phi}'$ for $\phi'$ and $\Phi$ for $\tilde{\Phi}$ (henceforth referred to as symmetric in tilde), so that
\begin{equation*}
\int_M \sum_{i,j} \frac{\partial \sigma_2}{\partial w_{ij}}(W) \; \phi' \; \mathrm{Hess}^{M,g}(\tilde{\Phi})(e_i, e_j) - \int_M \sum_{i,j} \frac{\partial \sigma_2}{\partial w_{ij}}(W) \; \tilde{\phi}' \; \mathrm{Hess}^{M,g}(\Phi)(e_i, e_j) = 0.
\end{equation*}
\label{IntegralSym}
\end{theorem}

\begin{proof}
From the definition of the Hessian, we will start from the integral
\begin{equation*}
\int_M \sum_{i,j} \frac{\partial \sigma_2}{\partial w_{ij}}(W) \; \tilde{\phi}' \; \Big( \langle \mathrm{grad}(\langle \mathrm{grad}(\Phi), e_j \rangle ), e_i \rangle - \langle \mathrm{grad}(\Phi), \nabla_{e_i} e_j \rangle \Big).
\end{equation*}
Now taking the $\frac{\partial \sigma_2}{\partial w_{ij}}(W) \; \tilde{\phi}'$ inside the inner products gives
\begin{equation*}
\int_M \sum_{i,j} \Big\langle \mathrm{grad}(\langle \mathrm{grad}(\Phi), e_j \rangle ), \frac{\partial \sigma_2}{\partial w_{ij}}(W) \; \tilde{\phi}' \; e_i \Big\rangle - \Big\langle \mathrm{grad}(\Phi), \frac{\partial \sigma_2}{\partial w_{ij}}(W) \; \tilde{\phi}' \; \nabla_{e_i} e_j \Big\rangle,
\end{equation*}
and integration by parts on the first term gives
\begin{equation*}
\int_M \sum_{i,j} -\langle \mathrm{grad}(\Phi), e_j \rangle \; \mathrm{div}\Big(\frac{\partial \sigma_2}{\partial w_{ij}}(W) \; \tilde{\phi}' \; e_i \Big) - \Big\langle \mathrm{grad}(\Phi), \frac{\partial \sigma_2}{\partial w_{ij}}(W) \; \tilde{\phi}' \; \nabla_{e_i} e_j \Big\rangle.
\end{equation*}
Now by the product rule we have 
\begin{multline*}
\int_M \sum_{i,j} -\langle \mathrm{grad}(\Phi), e_j \rangle \; \Big(\frac{\partial \sigma_2}{\partial w_{ij}}(W) \; \tilde{\phi}' \; \mathrm{div}(e_i) - \Big\langle \mathrm{grad}\Big( \frac{\partial \sigma_2}{\partial w_{ij}}(W) \; \tilde{\phi}'\Big), e_i \Big\rangle \Big) \\ - \Big\langle \mathrm{grad}(\Phi), \frac{\partial \sigma_2}{\partial w_{ij}}(W) \; \tilde{\phi}' \; \nabla_{e_i} e_j \Big\rangle.
\end{multline*}
Rearranging gives
\begin{multline}
\int_M \sum_{i,j} \langle \mathrm{grad}(\Phi), e_j \rangle \; \Big\langle \mathrm{grad}\Big( \frac{\partial \sigma_2}{\partial w_{ij}}(W) \; \tilde{\phi}'\Big), e_i \Big\rangle \\ - \Big\langle \mathrm{grad}(\Phi), \frac{\partial \sigma_2}{\partial w_{ij}}(W) \; \tilde{\phi}' \; \Big( \mathrm{div}(e_i) e_j + \nabla_{e_i} e_j \Big) \Big\rangle.
\label{Sym2}
\end{multline}

We will now examine the two terms in this integral separately. Starting with the first term, since $\mathrm{grad}$ satisfies Leibniz's rule, we have
\begin{equation*}
\int_M \sum_{i,j} -\langle \mathrm{grad}(\Phi), e_j \rangle \; \Big\langle \mathrm{grad}\Big( \frac{\partial \sigma_2}{\partial w_{ij}}(W) \Big) \tilde{\phi}' + \frac{\partial \sigma_2}{\partial w_{ij}}(W) \; \mathrm{grad}( \tilde{\phi}'), e_i \Big\rangle.
\end{equation*}
The linearity of the inner product and splitting the integral gives
\begin{multline*}
\int_M \sum_{i,j} -\langle \mathrm{grad}(\Phi), e_j \rangle \; \Big\langle \mathrm{grad}\Big( \frac{\partial \sigma_2}{\partial w_{ij}}(W) \Big) \tilde{\phi}', e_i \Big\rangle \\ + \int_M \sum_{i,j} -\langle \mathrm{grad}(\Phi), e_j \rangle \; \Big\langle \frac{\partial \sigma_2}{\partial w_{ij}}(W) \; \mathrm{grad}( \tilde{\phi}'), e_i \Big\rangle.
\end{multline*}
Factoring out $-\langle \mathrm{grad}(\Phi), e_j \rangle$ and $\tilde{\phi}'$ from the sum over $i$ in the first integral gives
\begin{multline}
\int_M \sum_{j} - \Big( \tilde{\phi}' \; \langle \mathrm{grad}(\Phi), e_j \rangle \; \sum_i \Big\langle \mathrm{grad}\Big( \frac{\partial \sigma_2}{\partial w_{ij}}(W) \Big), e_i \Big\rangle \Big) \\ + \int_M \sum_{i,j} -\langle \mathrm{grad}(\Phi), e_j \rangle \; \Big\langle \frac{\partial \sigma_2}{\partial w_{ij}}(W) \; \mathrm{grad}( \tilde{\phi}'), e_i \Big\rangle,
\label{Sym3}
\end{multline}
but then by Proposition \ref{gradSigma2} the first integral is zero. After factoring out $\frac{\partial \sigma_2}{\partial w_{ij}}(W)$ in the second integral, since $\tilde{\phi}'=\tilde{\Phi}$ we see that (\ref{Sym3}) is equal to
\begin{equation}
\int_M \sum_{i,j} \frac{\partial \sigma_2}{\partial w_{ij}}(W) \langle \mathrm{grad}(\Phi), e_j \rangle \langle \mathrm{grad}(\tilde{\Phi}), e_i \rangle.
\label{Sym4}
\end{equation}

Now we will look at the second term of (\ref{Sym2}). First let
\begin{equation*}
X=\sum_{i,j} \frac{\partial \sigma_2}{\partial w_{ij}}(W) \; \Big( \mathrm{div}(e_i) e_j + \nabla_{e_i} e_j \Big),
\end{equation*}
then the second term of (\ref{Sym2}) is equal to
\begin{equation*}
-\int_M \langle \mathrm{grad}(\Phi), \tilde{\phi}' X \rangle.
\end{equation*}
Trivially this is
\begin{equation*}
-\int_M \frac{1}{2}\langle \mathrm{grad}(\Phi), \tilde{\phi}' X \rangle + \frac{1}{2}\langle \mathrm{grad}(\Phi), \tilde{\phi}' X \rangle.
\end{equation*}
Then integration by parts on the second term gives
\begin{equation*}
-\int_M \frac{1}{2}\langle \mathrm{grad}(\Phi), \tilde{\phi}' X \rangle - \frac{1}{2}\Phi \mathrm{div}(\tilde{\phi}' X).
\end{equation*}
By the product rule, this is
\begin{equation*}
-\int_M \frac{1}{2}\langle \mathrm{grad}(\Phi), \tilde{\phi}' X \rangle - \frac{1}{2}\Phi \Big( \tilde{\phi}' \mathrm{div}(X) - \langle \mathrm{grad}(\tilde{\phi}'), X \rangle \Big).
\end{equation*}
Recall we have $\tilde{\phi}'=\tilde{\Phi}$, so after expanding we have
\begin{equation}
-\int_M \frac{1}{2}\langle \mathrm{grad}(\Phi), \tilde{\Phi} X \rangle + \frac{1}{2} \langle \mathrm{grad}(\tilde{\Phi}), \Phi X \rangle- \frac{1}{2}\Phi \tilde{\Phi} \mathrm{div}(X).
\label{Sym5}
\end{equation}

Plugging (\ref{Sym4}) and (\ref{Sym5}) into (\ref{Sym2}), we obtain
\begin{multline}
\int_M \sum_{i,j} \frac{\partial \sigma_2}{\partial w_{ij}}(W) \; \tilde{\phi}' \; \mathrm{Hess}^{M,g}(\Phi)(e_i, e_j) \\
\begin{aligned}
= \int_M \sum_{i,j} \Big[ & \frac{\partial \sigma_2}{\partial w_{ij}}(W) \langle \mathrm{grad}(\Phi), e_j \rangle \langle \mathrm{grad}(\tilde{\Phi}), e_i \rangle \\
&- \frac{1}{2} \Big\langle \mathrm{grad}(\Phi), \tilde{\Phi} \frac{\partial \sigma_2}{\partial w_{ij}}(W) \; \Big( \mathrm{div}(e_i) e_j + \nabla_{e_i} e_j \Big) \Big\rangle \\
&- \frac{1}{2} \Big\langle \mathrm{grad}(\tilde{\Phi}), \Phi \frac{\partial \sigma_2}{\partial w_{ij}}(W) \; \Big( \mathrm{div}(e_i) e_j + \nabla_{e_i} e_j \Big) \Big\rangle \\
&+ \frac{1}{2}\Phi \tilde{\Phi} \mathrm{div}\Big( \frac{\partial \sigma_2}{\partial w_{ij}}(W) \; \Big( \mathrm{div}(e_i) e_j + \nabla_{e_i} e_j \Big) \Big) \Big]
\label{Sym6}
\end{aligned}
\end{multline}

Note that since the $\{e_i\}$ are orthonormal, $W=w_{ij}$ is symmetric in $i$ and $j$. Now from the explicit formula for $\frac{\partial \sigma_2}{\partial w_{ij}}(W)$ in (\ref{Sigma2Deriv}), we have that the first term on the right hand side of (\ref{Sym6}) is symmetric in tilde. Clearly the middle two terms together are symmetric in tilde and the last term is as well. Therefore the whole integral is symmetric in tilde.
\end{proof}

\section{The Rigidity Theorem}

Denote by $dS^{n+1, +}$ the region of de Sitter space for which $\rho$ is positive. We can now prove the main rigidity result following the idea for the Riemannian case in \cite{GuanShen}.

\begin{theorem}
Let $M$ and $\tilde{M}$ be two spacelike hypersurfaces in $dS^{n+1,+}$ with $K>\bar{K}$ such that there exists a local isometry $f:M\to \tilde{M}$. Then $f$ is the restriction of some global isometry $F$ of de Sitter space.
\end{theorem}

\begin{proof}
We start with the integral equations in Theorem \ref{IntegralEqn}, subtracting (\ref{Integral0a}) from (\ref{Integral0c}) and (\ref{Integral0b}) from (\ref{Integral0d}), by Theorem \ref{IntegralSym} the right hand sides will cancel and we obtain
\begin{subequations}
\begin{equation}
\int_M \tilde{\phi}' \sigma_2(W) \langle V, \nu \rangle = \int_M \phi' \sigma_{1,1}(\tilde{W},W) \langle \tilde{V}, \tilde{\nu} \rangle,
\label{Rigidity1a}
\end{equation}
\begin{equation}
\int_M \phi' \sigma_2(\tilde{W}) \langle \tilde{V}, \tilde{\nu} \rangle = \int_M \tilde{\phi}' \sigma_{1,1}(W,\tilde{W}) \langle V, \nu \rangle.
\label{Rigidity1b}
\end{equation}
\end{subequations}
Now note that $\sigma_2(W)=\sigma_2(\tilde{W})$ and $\sigma_{1,1}(W,\tilde{W})=\sigma_{1,1}(\tilde{W},W)$, so adding (\ref{Rigidity1a}) and (\ref{Rigidity1b}) and rearranging we obtain
\begin{equation*}
\int_M \Big( \tilde{\phi}' \langle V, \nu \rangle + \phi' \langle \tilde{V}, \tilde{\nu} \rangle \Big) \Big( \sigma_2(W) - \sigma_{1,1}(W,\tilde{W}) \Big) = 0.
\label{Rigidity2}
\end{equation*}
Since $M, \tilde{M}\subset dS^{n+1,+}$, $\phi'$ and $\tilde{\phi}'$ are positive on all of $M$. Furthermore since $M$ and $\tilde{M}$ are spacelike, $\langle V, \nu \rangle$ and $\langle \tilde{V}, \tilde{\nu} \rangle$ are strictly less than zero. So the quantity $\tilde{\phi}' \langle V, \nu \rangle + \phi' \langle \tilde{V}, \tilde{\nu} \rangle$ is strictly less than zero. Now from the inequality in Theorem \ref{GardingIneq} we have that $\sigma_2(W) - \sigma_{1,1}(W,\tilde{W})$ is less than or equal to zero. Therefore the only way the integral in (\ref{Rigidity2}) can be zero is if $\sigma_2(W) - \sigma_{1,1}(W,\tilde{W})$ is identically zero on all of $M$.

Now from the equality result of G\r arding's inequality for hyperbolic polynomials in \S\ref{HyperbolicRigidity}, since $\sigma_2(W) - \sigma_{1,1}(W,\tilde{W})=0$ is zero on all of $M$, we have that $W=\tilde{W}$. Since the first and second fundamental forms are preserved under the map $f$, it must be the restriction of some global isometry $F$, hence the theorem is proved.
\end{proof}

\bibliographystyle{amsplain}
\bibliography{PhDBib}{}

\end{document}